\documentclass{paper}

\usepackage{arxiv}

\usepackage{amsmath,amsfonts}
\usepackage{array}
\usepackage[caption=false,font=normalsize,labelfont=sf,textfont=sf]{subfig}
\usepackage{textcomp}
\usepackage{stfloats}
\usepackage{verbatim}
\usepackage{graphicx}
\usepackage{cite}
\usepackage{float}
\usepackage[dvipsnames,table]{xcolor}
\usepackage[pagebackref,breaklinks,colorlinks]{hyperref}
\usepackage{url}

\definecolor{color_todo}{rgb}{.6,.8,.2}

\newcommand{\email}[1]{\nolinkurl{#1}}

\newcommand{\NN}{\mathbb{N}}
\renewcommand{\AA}{\mathcal{A}}
\newcommand{\R}{\mathcal{R}}

\definecolor{tablelighterblue}{HTML}{DAE8FC}
\definecolor{tabledarkerblue}{HTML}{aacaf8}

\usepackage[capitalize]{cleveref}

\usepackage{multirow}
\usepackage{amssymb}
\usepackage{algcompatible}
\usepackage{algpseudocode}
\usepackage{algorithm2e, setspace}

\usepackage{caption}

\title{Unrolled Iterative Reconstruction Networks
for Implicit Motion Compensation in
Dynamic Nano Computed Tomography}

\author{
  Alice Oberacker \\
  Department of Mathematics\\
  Saarland University\\
  Saarbr\"ucken, Germany \\
  \texttt{alice.oberacker@num.uni-sb.de} \\
       \And
  Anne Wald \\
  Institute for Numerical and Applied Mathematics \\
  University of G\"ottingen 
  G\"ottingen, Germany \\
  \texttt{a.wald@math.uni-goettingen.de} \\
       \And
  Thomas Schuster \thanks{Corresponding Author}\\
  Department of Mathematics\\
  Saarland University\\
  Saarbr\"ucken, Germany \\
  \texttt{thomas.schuster@num.uni-sb.de} \\
}

\begin{document}

\maketitle

\begin{abstract}
Nano-CT is an essential imaging tool that is highly susceptible to motion-induced artifacts caused by vibrations and thermal drifts. These dynamic perturbations introduce significant model inexactness into the reconstruction process. This paper proposes and evaluates two unrolled neural network architectures designed to learn the reconstruction mapping directly from perturbed measurement data by incorporating the mathematical framework of the Sequential Subspace Optimization (SESOP) algorithm. The Fully Learned Operator (FLO) model approximates the true inverse mapping from only the perturbed measurement data, while the Radon network serves as a physics-informed model that explicitly incorporates the forward and adjoint Radon operators within its unrolled iterations.
Through experiments involving simulated vibrations and linear shifts, these learned models were compared against the Filtered Back Projection, the Dremel method, and the SESOP-Kaczmarz technique. The results demonstrate that the unrolled architectures significantly outperform traditional methods in overall performance and noise suppression. Specifically, the FLO model proved to be the most robust, successfully recovering internal features and recognizable object shapes even under extreme linear shifts where non-learned methods failed. This model provided superior visual clarity and accurately reconstructed densities, establishing it as a superior high-performance solution for dynamic CT imaging in environments with significant modeling inexactness.

\end{abstract}



\section{Introduction}

Computed tomography (CT) is a well-established medical imaging technique that also finds application in the industrial setting. CT imaging at the nano scale (nano-CT) is used to visualize the inner structures of small objects. For example, in material testing, nano-CT is particularly important as a resolution of up to 100 nanometers is used for the development of microchips \cite{circuit_testing, cryst11060677}; electrode degradation studies in batteries \cite{cryst11060677} and many other applications. \\

Due to the small size of the scanned objects, noise has an amplified impact on the final imaging quality.
Unfortunately, these small objects are highly susceptible to ambient conditions. For example, the CT scanner can cause vibrations and temperature changes leading to \emph{thermal drifts}. Both of these motions cause a misalignment between the object and the X-ray projections \cite{Liu, oliveira2014}, which degrades the measurement data and causes motion blur in the reconstructions when classical reconstruction methods are applied.

To mitigate the effects of object movement, researchers have proposed post-processing solutions that estimate motion from the data itself \cite{Huang, hahn2017motion}; active monitoring techniques that record displacement throughout the scan \cite{Chen, Burger_2017}. 
Reconstructing the cross-sectional view of an object using CT measurement data requires solving an inverse problem. 
Since the object is scanned from multiple angles sequentially, the motion of the object causes the object to be in different states for every projection angle. Thus, the resulting measurement data are inconsistent with any individual state of the object \cite{hahn2014reconstruction, hahn2017motion, Blanke_2020}.

In this article, we consider machine learning-based inversion methods. More specifically, we propose two different architectures, both inspired by unrolled iterative regularization schemes. One is the fully learned operator (FLO) network, which is trained from solution-data pairs with no additional information on the forward operator. The second is the Radon network, which replaces the unknown (motion-dependent) forward mapping in the unrolled scheme by the Radon transform, which is seen as an appropriate approximation. Numerical experiments prove that the proposed unrolled networks outperform existing methods as the classical Filtered Backrojection (FBP) scheme or the Dremel method.

\section{Problem definition and setting}

\noindent Nano-CT refers to the tomographic imaging of an object's density in the nanometer and lower micrometer range, a field that has pushed detail resolution to 100 nm \cite{LUTTER202110}. The object is illuminated by an X-ray beam with intensity $I_0$, which undergoes attenuation based on the absorption coefficient $f: \mathbb{R}^2 \rightarrow \mathbb{R}$ of the medium \cite{Natterer_img_recon}. Following the Beer-Lambert law, the decrease in intensity is proportional to the traveled distance and the density of the material, leading to the relation
\begin{equation*}
    \ln \left(\frac{I_0}{I_1}\right) = \int_{L_{\omega(\varphi), s}} f(x)\, \mathrm{d}x
\end{equation*}
between the measured intensities $I_0, I_1$ and the unknown target $f$. The beam path $L_{\omega(\varphi), s}$ is a line parameterized by the detector position $s\in\mathbb{R}$ and the normal vector $\omega(\varphi) := (\cos \varphi, \sin \varphi)^T$, where $\varphi$ is the scan angle \cite{Natterer, Natterer_img_recon}, i.e., $L_{\omega(\varphi), s} = \{s\omega (\varphi)+t\omega (\varphi)^\perp:t\in\mathbb{R} \}$. The mapping
\begin{equation} \label{eq:radon_def}
    (\R f)(\omega, s) := \int_{L_{\omega(\varphi), s}} f(x) \, \mathrm{d}x
\end{equation}
is the Radon transform, providing the mathematical model of 2D-CT \cite{Natterer_img_recon}. 
Conventional scanning geometries in two dimensions are the parallel beam or fan beam setting. Many nano-CT systems utilize a fan beam geometry where the source emits a fan of radiation \cite{Graetz_2021}. The underlying principle represents an inverse problem that consists of reconstructing the cross-sectional phantom from the collected sinogram data. The extreme magnification required for nanoscale resolution renders the reconstruction process highly sensitive to even small disturbances (ill-posedness). During the measurement process, which is not instantaneous, the assumption of a static object $f$ often fails. Instead, the object undergoes changes due to unavoidable relative motion between the sample and the tomograph, caused by mechanical vibrations, manufacturing tolerances in the rotating stage, or thermal drift caused by the heat-induced expansion of scanner components \cite{oliveira2014, Liu, Dremel_2018}. Mathematically, this implies that the object $f$ depends not only on the spatial variable $x$ but also on time $t$. Failure to compensate for this motion results in severe artifacts, as the time-dependent measurements $g(t)$ become inconsistent with the static operator $\R$ \cite{Burger2023IllposednessOT}. To this end the state of the object $f$ at any time $t$ is mathematically described using a \emph{deformation model}.
\subsection{The Deformation Model}

\noindent Let $\mathcal{B} \subset \mathbb{R}^2$ be a bounded domain representing the field of view. Following Blanke et al. \cite{Blanke_2020}, we describe the temporal evolution of the object through a deformation model $\Gamma : [0, T] \times \mathcal{B} \rightarrow \mathcal{B}$.
For a fixed time $t$, $\Gamma(t, \cdot)$ is a spatial transformation that maps a point $x$ in the reference state (usually $t=0$) to its position $\hat{x}$ in the deformed state. The time-dependent object is thus represented as the composition
\begin{equation}
    f_t(x) = (f \circ \Gamma)(t, x) = f(\Gamma(t, x)) = f(\hat{x}).
\end{equation}

\subsection{Construction of the Dynamic Operator}

\noindent The exact (dynamic) forward operator $\mathcal{A}_\Gamma$ incorporates this motion directly into the physics of the Radon transform. The dynamic process is then given as
\begin{align}
    (\mathcal{A}_\Gamma f)(t) &= (\R f_t)(t) \nonumber\\
    &= \int_{L(\omega_t, s)} f(\Gamma(t, x)) \, dx \\
    &= g(t),  
\end{align}
where $\omega_t = \omega(\varphi(t))$ denotes the source depending on time $t$.
Here, each projection $g(t)$ corresponds to a unique angle $\varphi(t)$, but is applied to a different state of the deformed object $f_t$. In this sense, the acquired data $g(t)$ are inconsistent, meaning that $g(t)\not\in \mathrm{Image} (\R)$.

\subsection{Operator Inexactness and Modeling Errors}

\noindent In practice, $\Gamma$ and, thus, $\mathcal{A}_\Gamma$ are not given. Hence, we follow the lines in \cite{Blanke_2020} and approximate $\mathcal{A}_\Gamma$ by a mapping $\mathcal{A}^\eta$, i.e., $\mathcal{A}^\eta \approx \mathcal{A}_\Gamma$, where $\eta>0$ denotes the model inaccuracy,
$$\|\mathcal{A}_\Gamma - \mathcal{A}^\eta\| \leq \eta.$$
For a bounded object $\|f\| \leq \rho$, the total modeling error is then bounded by $\eta \rho$. Accordingly, we denote the measurement data by $g^\eta$ and the inverse problem to be solved is $\mathcal{A}^\eta f=g^\eta$.

\section{Materials and Methods}

\noindent Applying conventional inversion techniques like Filtered Back Projection (FBP) \cite{Natterer} to dynamic measurements results in significant motion artifacts. These are particularly detrimental in nano-CT, where thermal drift or mechanical vibrations \cite{oliveira2014, Liu} cause misalignments between projections and the object. Several strategies have been proposed to mitigate these effects:
\begin{itemize}
    \item Motion Estimation: Estimating the movement directly from the sinogram data \cite{Katsevich, hahn2017motion, Ramos17, Huang}.
    \item Direct Measurement: Using markers or sensors to track motion during the scan \cite{Burger_2017, Chen}.
    \item Dynamic Reconstructions: Reconstructing a single state of the dynamic object (e.g., $t=0$) from all inconsistent data $g(t)$, as explored by Burger et al. \cite{Burger2023IllposednessOT}.
\end{itemize}
We pursue another approach by using neural networks that are trained by supervised learning and whose architecture is inspired by unrolled iterative regularization methods. This way we combine stability with accuracy in the reconstruction process.

\subsection{Sequential reconstruction methods}

\subsubsection{SESOP } 
We consider sequential subspace optimization \cite{SESOP1, SESOP2}, which can be seen as a generalization and acceleration of the Landweber method \cite{landweber1951}.
The Landweber method updates the current iterate $f_n$ using the negative gradient of he residuum as search direction. The SESOP method extends this idea by using more than one gradient simultaneously in each iteration, which improves convergence.
A SESOP iteration step applied to $\mathcal{A}^\eta f=g^\eta$ reads as  
\begin{equation}
    f_{n+1} = f_n - \sum_{i \in I_n} t_{n,i} (\AA^\eta)^* (\AA^\eta f_i - g^\eta)),
        \label{eq:sesop_method}
\end{equation}
with $f_n$ denoting the current iterate, $I_n\subset\{0,\ldots,n\}$ is a finite index set, $t_{n,i}>0$ is the step size, and $(\AA^\eta)^*$ is the adjoint operator to $\AA^\eta$ (with respect to appropriate Hilbert spaces $X$, $Y$).

\subsubsection{The SESOP-Kaczmarz method} 
The SESOP-Kaczmarz method is a hybrid reconstruction scheme that merges the acceleration capabilities of SESOP~\cite{SESOP1, SESOP2} with the structural flexibility of the Kaczmarz method~\cite{Natterer}. While standard SESOP speeds up convergence by utilizing multiple search directions simultaneously, the Kaczmarz-method is capable to handle semi-discrete forward mappings. By integrating the Kaczmarz approach, the inverse problem can be decomposed into sub-problems defined by specific scan angles $\varphi$ and detector points $s$. This decomposition is particularly advantageous for dynamic CT, where modeling errors are not static but vary over time and, thus, rotation angle $\varphi$.

\noindent To iterate through the sub-problems, two index sets are introduced. The index set for the scan angles $\varphi_k$ is given by $\mathcal{K} := \{0, ..., K-1 \}$ where $\varphi_k \in [0, \pi]$ for $k \in \mathcal{K}$; similarly, the index set for the detector points $s_l$ is defined as $\mathcal{L} := \{0, ..., L-1 \}$ with $s_l \in [-1,1]$ for $l \in \mathcal{L}$. 
The algorithm updates the current iterate $f_n$ as
$$f_{n+1} = f_n - \sum_{i \in I_n} t_{n,i} (\AA_{[i]}^\eta)^* (\AA_{[i]}^\eta f_n- g^\eta_{[i]}),$$ 
with $[i] = i \mod K \cdot L$.\\ 
This semi-discrete realization of the operator allows the method to incorporate local characteristics of the inverse problem, such as local quantum noise $\delta_j$ and local model inexactness $\eta_j$. By separating these errors into sub-problems, the reconstruction process remains robust even when the dynamic forward operator deviates from the ideal model.\\
Regularization within the SESOP-Kaczmarz framework is achieved by metric projections onto stripes that contain the set of solutions of $\AA^\eta f=g^\eta$. Moreover, these stripes
represent a feasible manifold that accounts for the combined influence of measurement noise and model error, see \cite{SESOP2, Blanke_2020} for details. The iteration process is stopped with the smallest $n$ such that the discrepancy principle
$$\|\mathcal{A}^\eta_{[i]} f_n - g^\delta_{[i]} \| \leq \tau(\eta_{[i]} \rho + \delta_{[i]})$$
is satisfied where $\tau>1$ acts as a tolerance parameter. This ensures that the algorithm terminates once the reconstruction is consistent with the known noise level and modeling bounds, preventing the over-amplification of artifacts that often plagues non-regularized iterative methods in dynamic settings. For a detailed analysis of this method we refer to \cite{Blanke_2020} and for implementation details and settings to \cite{Oberacker_2025}.

\subsubsection{The Dremel method} 
The Dremel method, introduced by Kilian Dremel in 2018 \cite{Dremel_2018}, is a model-based iterative reconstruction procedure designed to address the challenges of rigid body motion in CT. This algorithm is fundamentally rooted in the Kaczmarz iteration, but it distinguishes itself by interlacing a specialized correction step that dynamically adjusts the semi-discrete forward operator $\mathcal{A}_k^\eta$ and its adjoint during the reconstruction process. This mechanism is primarily intended to identify and compensate for shifts within the detector plane, though it is also capable of estimating movements of the object or the X-ray source by translating these physical shifts into the operator's geometry. During each main iteration, the algorithm cycles through all scan angles and optimizes the alignment between the measured data and the forward projection of the current intermediate reconstruction. This optimization is achieved through a cross-correlation process that determines the ideal detector shift required for the next iteration step. 
While the Dremel method has shown significant empirical success in handling modeling errors and dynamic inconsistencies \cite{NN_inexactOp}, it is important to note that it currently lacks a formal mathematical convergence analysis, unlike more established iterative frameworks such as SESOP or the standard Kaczmarz method. The pseudocode for this method and the method configurations can be found in~\cite{NN_inexactOp, Oberacker_2025}.

\subsection{Unrolled neural network}

\begin{figure}[h]
  \centering
  \includegraphics[width=.8\linewidth]{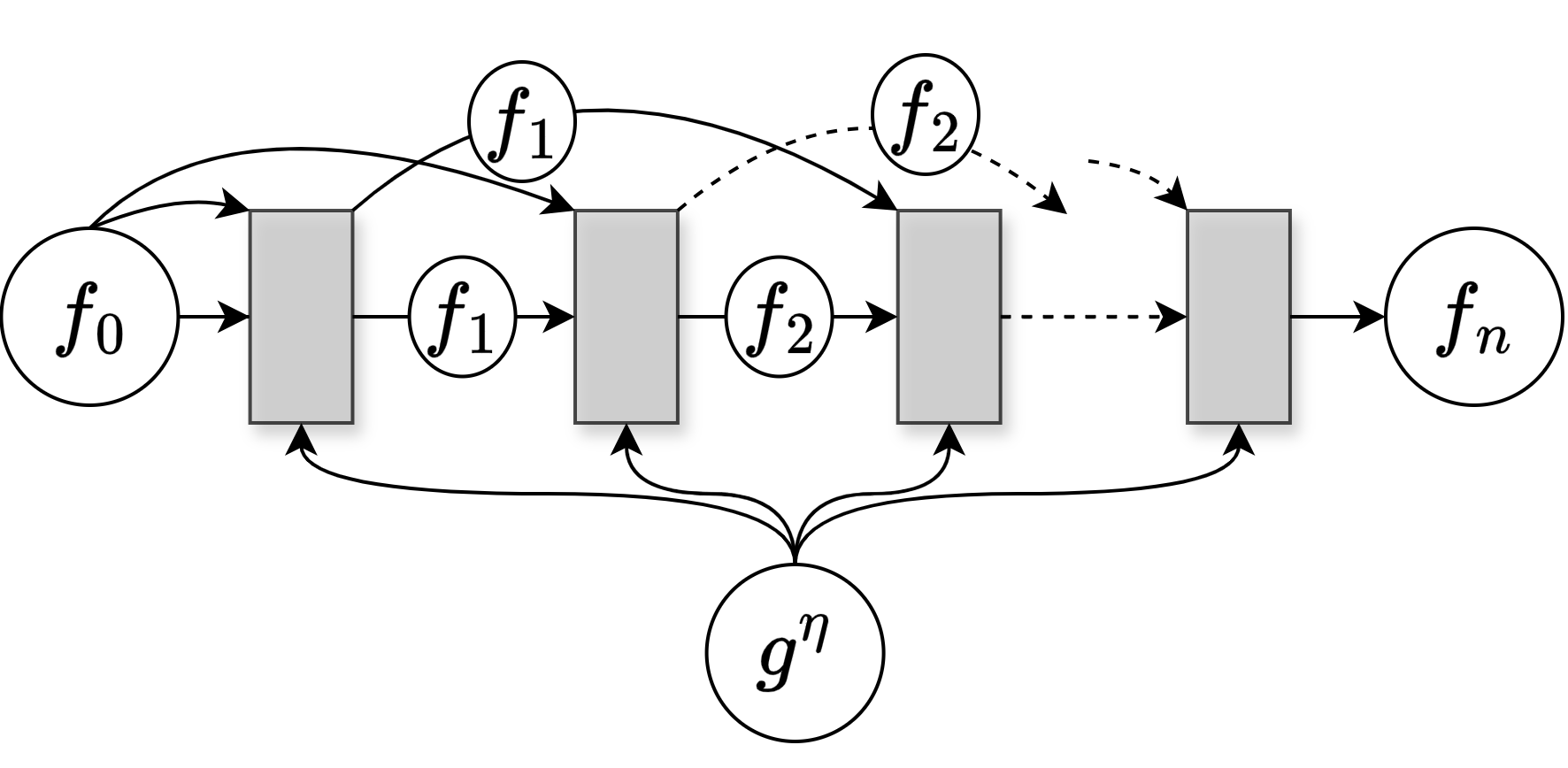}
  \caption{Overall structure of SESOP iteration for two search directions and $n$ unrolled iterations.}
\label{fig:sesop_overview_architecture}
\end{figure}

We choose a data-driven approach relying on a neural network architectures capable of learning the reconstruction mapping $f$ directly from the perturbed data $g^\eta$ without explicit prior knowledge of~$\Gamma$. 
Many neural network architectures have been successfully applied to different tasks in the field of CT \cite{Baguer_2020, Baguer2020ComputedTR, CNNreconstruction}. A popular approach consists of applying a non-learned reconstruction method followed by a post-processing network~\cite{NN_inexactOp}. Unrolled neural networks \cite{GregorL10, monga} encapsulate the reconstruction process by implementing the iterative steps of a sequential method. Let $\Phi_n: X \times Y \rightarrow X$ be a functional defining an arbitrary iterative scheme to approximate $f$ from $g^\eta$, that means $f_{n+1} = \Phi_n(f_n,g^\eta)$. For a fixed~$N \in \NN$ the $N$-th iterate is calculated by sequentially applying~$N-1$ iterations to an initial guess $f_0$. This can be written as $f_N = \Phi_{N-1} \circ \Phi_{N-2} \circ \dots \circ \Phi_0(f_0, g)$. Each iteration step~$\Phi_n$ is now approximated by a set of neural network layers~$\Psi_n\approx \Phi_n$. For a fixed number $N$ of iterations the unrolled neural network then reads as 
\begin{equation}\label{eq:unrolled_template}
  f_N \approx \Psi_{N-1} \circ \Psi_{N-2} \circ \dots \circ \Psi_0(f_0, g).
\end{equation}
The advantage of such networks is that they encapsulate the mathematical framework of iterative methods while leveraging the adaptive power of learned methods. Once trained, the inference or reconstruction time is faster compared to the given iteration scheme $\Phi_n$.
A key idea of applying such a method is that a neural network can learn from perturbed data $g^\eta$ and compensate for the motion $\Gamma$. By approximating the SESOP iteration, the network retains some structure from the sequential method while learning network parameters from the measurement data.
The network therefore integrates the inexactness of the inverse problem into its reconstruction process.
Combining the iteration step of SESOP shown in \eqref{eq:sesop_method}, with \eqref{eq:unrolled_template}, the network block thus reads as
\begin{equation}\label{eq:network_block}
    \Psi_n((f_i)_{i \in I_n}, g^\eta) = f_n - \sum_{i \in I_n} t_{n,i} \AA^* (\AA f_i - g^\eta)
\end{equation}
with $\AA$ now denoting a linear mapping defining the trainable network parameters. Note that the number of inputs for each network block depends on the index set $I_n$. The general structure of these iteration steps can be seen in \cref{fig:sesop_overview_architecture}.

\section{Architectures}

Two architectures are presented, both constructed by unrolling the SESOP algorithm: 1.) The \emph{fully learned operator network} (FLO network) learns the reconstruction provided only with the (motion) perturbed measurement data $g^\eta$. That means, the FLO network only learns from training data with no a priori information about the unknown (motion-dependent) operator $\AA_\Gamma$. 2.) The \emph{Radon network}, on the other hand, additionally receives the (adjoint) Radon transform, i.e., $\AA$ in \eqref{eq:network_block} is replaced by $\R$. Both networks additionally use an initial reconstruction $f_0$ which is computed by the FBP. Architecture and training details can be found in Supplementary Material~\ref{app:training_details} and~\ref{app:arch_details}.

\subsection{Fully Learned Operator Network}

\noindent The FLO Network \cite{Oberacker_2025} is designed to approximate a solution of $(\AA_\Gamma f)(t) = g(t)$
by simulating the SESOP iteration through a series of dedicated neural network blocks. It is, therefore, referred to as $\mathcal{N}_\theta^\AA$. This architecture is derived directly from the SESOP iteration step. 
To explain the structure of the network, the SESOP algorithm is discussed for two search directions, $I_n = \{n-1,n\}$. Then one iteration step reads as follows:
\begin{align}
        f_{n+1} &= f_n &-&  t_{n,n} \AA^* (\AA f_n - g^{\eta}) -  t_{n,n-1} \AA ^* (\AA f_{n-1} - g^{\eta}) \nonumber \\
                &= f_n &-& t_{n,n} \AA^* \AA f_n \nonumber \\
                  & &-& t_{n,n-1} \AA^* \AA  f_{n-1}  \nonumber \\
                  & &+& (t_{n,n} + t_{n,n-1}) \AA^* g^{\eta}.
                \label{eq:resesop_nn}
\end{align}
The individual components of the iteration step are approximated by collections of network layers. The iterate $f_{n+1}$ can then be approximated
\begin{equation}
    f_{n+1} \approx f_n - \mathcal{N}^{\text{sd}}_{\theta^1_n} (f_n) - \mathcal{N}^{\text{sd}}_{\theta^2_n} (f_{n-1}) + \mathcal{N}^{\text{D}}_{\theta^0_n}(g^{\eta}),
\end{equation}
where neural networks $\mathcal{N}^{\text{sd}}_{\theta^1_n} (f_n)$, $\mathcal{N}^{\text{sd}}_{\theta^2_n}(f_{n-1})$ and $\mathcal{N}^{\text{D}}_{\theta^0_n}(g^{\eta})$ approximate $t_{n,n} \AA^* \AA f_n$, $t_{n,n-1} \AA^* \AA  f_{n-1}$, $(t_{n,n} + t_{n,n-1}) \AA^* g^{\eta}$, respectively.
To optimize computational efficiency, the network performs 
the backprojection $\mathcal{A}^* g^\eta$ only once in a static pre-processing block, $\mathcal{N}^{\text{d}}_{\theta^0}(g^{\eta})$, outside the main iteration loop, while lighter internal blocks, $\mathcal{N}^{f}_{\theta^0_n}$, refine these results during each unrolled step. Thus, 
$$(t_{n,n} + t_{n,n-1}) \AA^* g^{\eta} \approx \mathcal{N}^{\text{D}}_{\theta^0_n}(g^{\eta}) = \mathcal{N}^{f}_{\theta^0_n} \circ \mathcal{N}^{\text{d}}_{\theta^0}(g^{\eta}).$$ Note that $\mathcal{N}^{\text{d}}_{\theta^0}$ is independent of $n$ as these layers are applied outside of the iteration steps and the result is used for each $\mathcal{N}^{f}_{\theta^0_n}$ with $n \in I_n$. A visualization of this can be found in \cref{fig:sesop_block_architecture_v2}. More detailed illustrations of this architecture can be found in Supplementary Material~\ref{app:N_0(g)} and~\ref{app:detail_sesop_block} as well as \cite{Oberacker_2025}.
\begin{figure}[h]
  \centering
  \includegraphics[width=.99\linewidth]{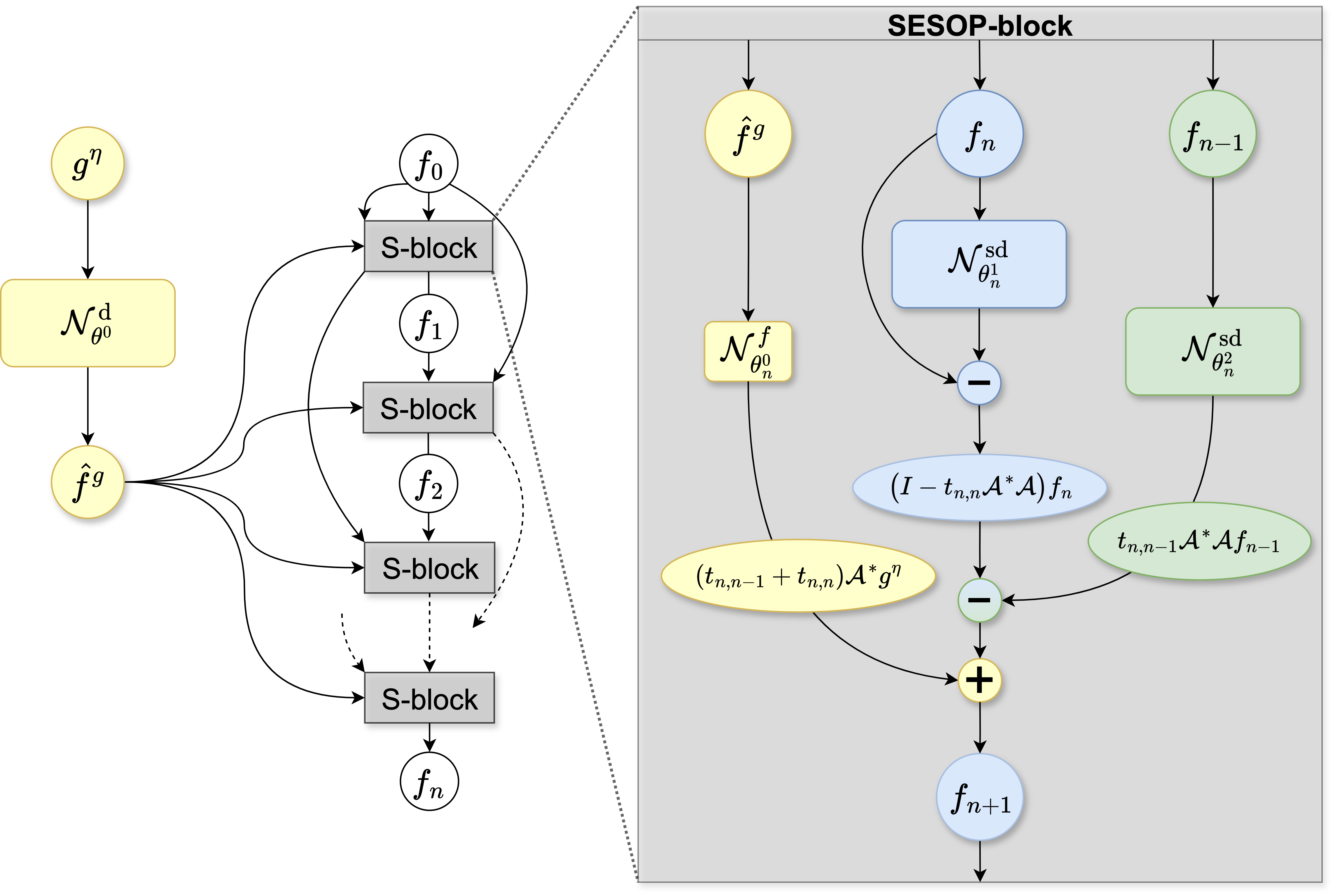}
  \caption{The SESOP-block structure includes three CNN blocks, one for each search direction and one single CNN layer for the measurement data, visualized by a smaller rectangle. Outside of the iteration block the main transformation for the measurement data is performed. }
  \label{fig:sesop_block_architecture_v2}
\end{figure}

The FLO Network is characterized by its capacity to approximate the inverse mapping $\mathcal{A}_\Gamma^{-1}$ entirely through data-driven blocks without requiring explicit knowledge of the forward operator. Its primary structural innovation lies in the decomposition of the SESOP iteration into independent sub-problems where the processing of measurement data $g^\eta$ is structurally isolated from the refinement of search directions $f_n$ and $f_{n-1}$. The architecture relies on deep residual blocks composed of convolution, batch normalization, and activation layers to mitigate gradient degradation and maintain high-frequency image fidelity across many unrolled iterations.

\subsection{The Radon Network}

\begin{figure}[h]
  \centering
  \includegraphics[width=0.99\linewidth]{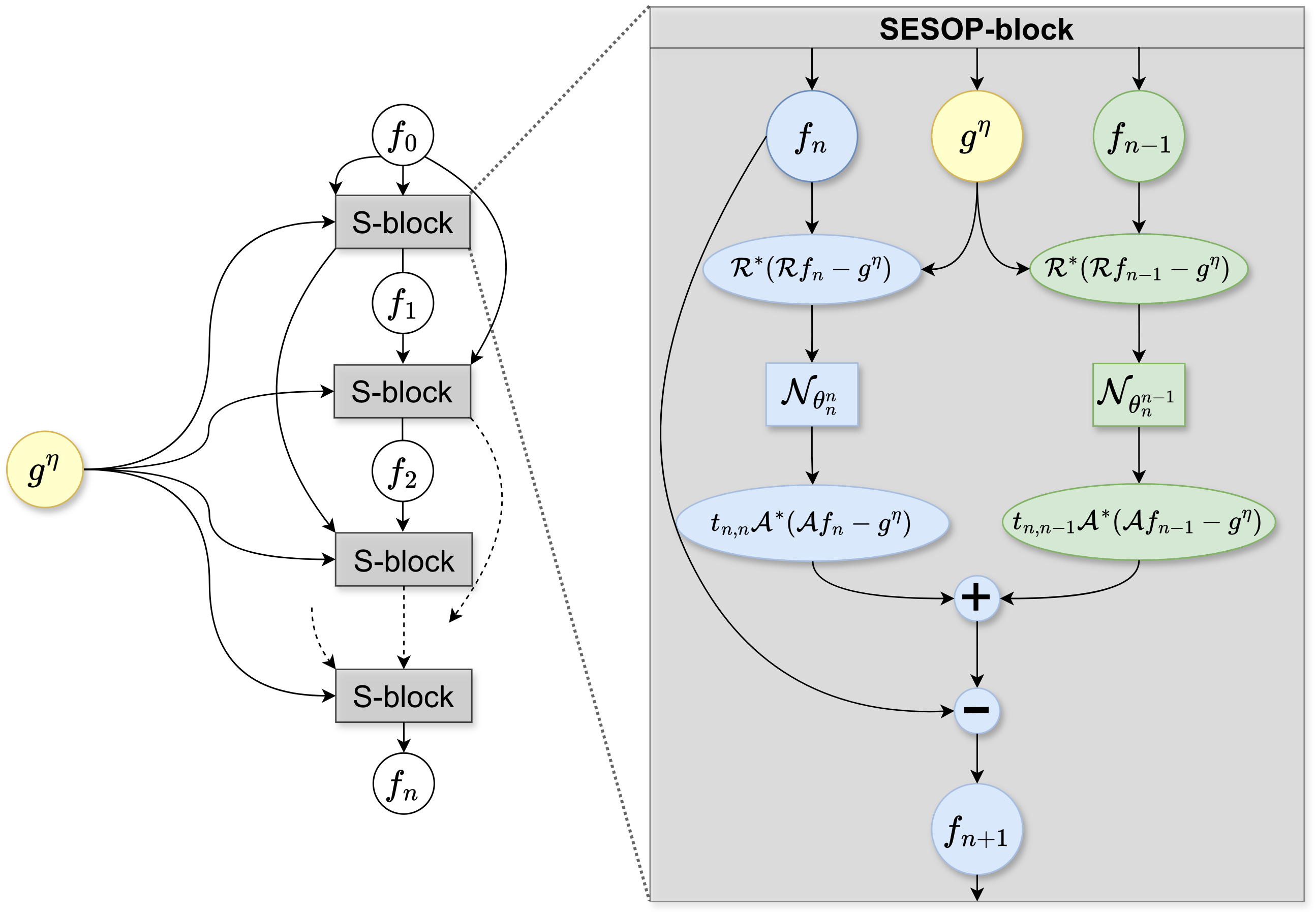}
  \caption{The structure of the SESOP-block includes one CNN block per search direction. The input of each of these block consists of the output of a previous iterations and the measurement data $g^\eta$. }
  \label{fig:adjoint_sesop_block}
\end{figure}

\noindent The Radon Network \cite{Oberacker_2025} serves as a physics-informed alternative to the FLO model by incorporating the Radon transform directly into the architecture. Therefore, this network is referred to as $\mathcal{N}_\theta^\R$. 
The SESOP iteration for the Radon equation is given by
\begin{equation}
    f_{n+1} = f_n - \sum_{i \in I_n} t_{n,i} \R^* (\R f_i - g^\eta)).
        \label{eq:sesop_method_radon}
\end{equation}
The iterate $f_{n+1}$ is then approximated through a set of network layers for each search direction
\begin{equation}
    f_{n+1} \approx f_n - \sum_{i \in I_n} \mathcal{N}_{\theta^i_n}(\mathcal{R}^*,\mathcal{R}, f_i, g^{\eta}),
\end{equation}
where $\mathcal{N}_{\theta^i_n}(\mathcal{R}^*,\mathcal{R}, f_i, g^{\eta})$ approximates $t_{n,i} \R^* (\R f_i - g^\eta))$ for each search direction with index $i$, the measurement data $g^\eta$ and the respective iterate $f_i$.
The structure of this architecture is visualized in \cref{fig:adjoint_sesop_block}.

The Radon Network serves as a physics-informed alternative that incorporates the Radon transform or its adjoint directly into the network's input pipeline. Its defining characteristic is the shift from learning the operator physics from the measurement data alone to learning artifact compensation, as it applies the adjoint Radon transform, a static operator, to the dynamic measurement data. By providing the network with the inexact operator, the architecture is primed with the general structural features of the object, allowing the convolutional layers to focus specifically on unlearning the streaks and blurs caused by object motion. This framework processes search directions in conjunction with the measurement data, effectively learning to map the inconsistencies of the inexact operator~$\mathcal{A}^\eta$ back to a motion-free ground truth while benefiting from the spatial weight-sharing efficiency of convolutional neural networks.

\subsection{Architecture Details}
\noindent Both network designs are based on similar deep learning techniques. All convolutional layers are paired with batch normalization (see \cite{batchnorm}) followed by activation layers. This combination of layers is called a CBA layer. By normalizing the activations of a mini-batch to a standard distribution and utilizing learnable parameters to maintain representational power, batch normalization ensures that the network remains robust throughout the optimization process. The leaky ReLU function was used for internal non-linearities with a negative slope of $0.01$. To ensure that the final reconstruction was within the valid physical range of $[0, 1]$, the sigmoid activation function was applied to the output layer.

\begin{figure}[H]
    \centering
    \includegraphics[width=3in]{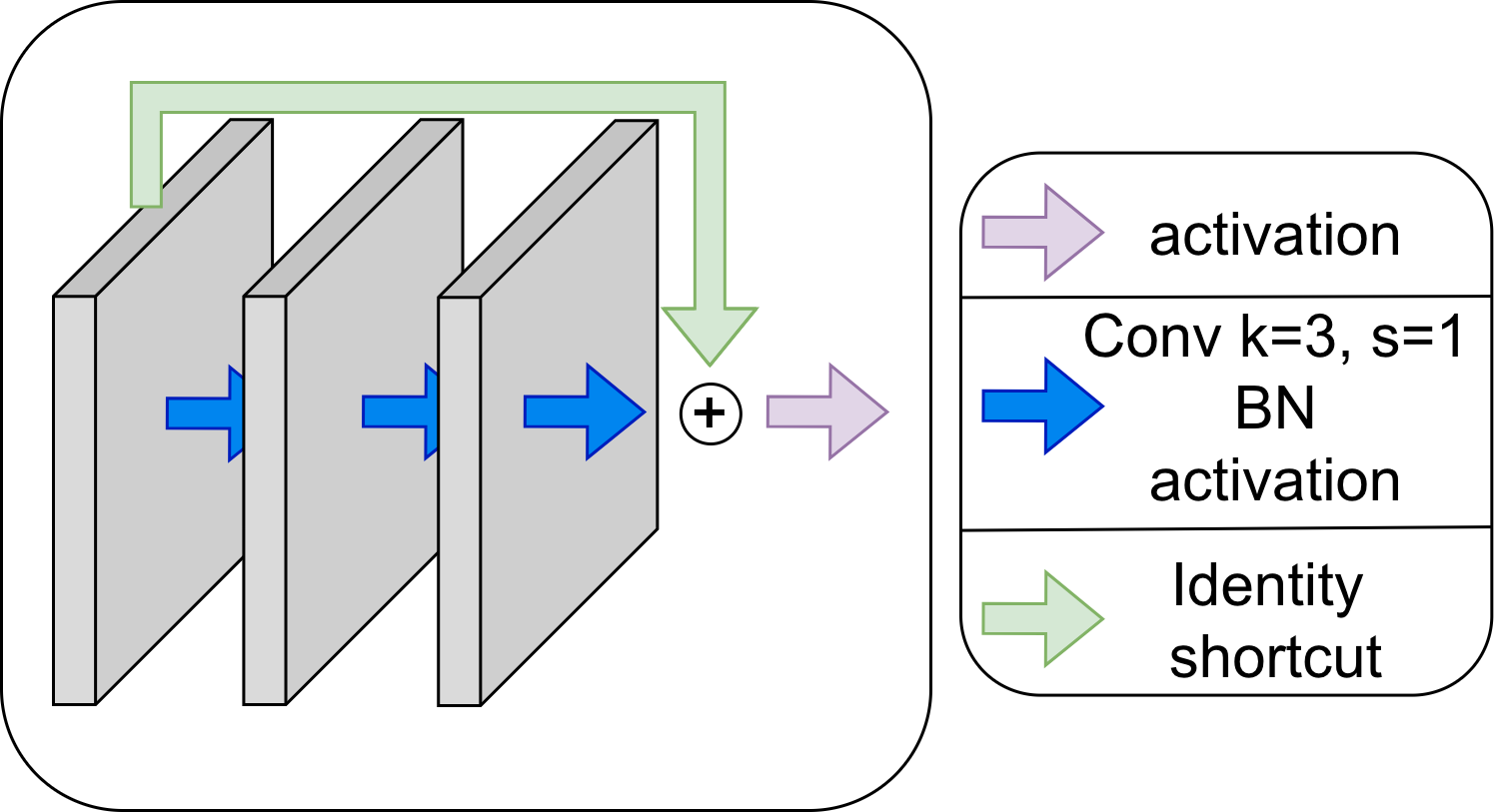}
    \caption{Implementation of residual layer. Three convolutional layers are combined with batch normalization and activations; the shortcut is implemented by PyTorch's identity layer. }
    \label{fig:res_layer}
\end{figure}

To address the degradation problem in deep networks, where increasing depth leads to higher training errors despite theoretically improved representational capacity \cite{deep_res} residual learning blocks were used. Note that in our setting depth equals the number of unrolled iteration steps. While deeper architectures are often necessary to learn complex features \cite{very_deep_cnn, deep_cnn}, standard non-linear layers struggle to approximate simple identity mappings during optimization \cite{deep_res}. The deep residual learning framework overcomes this by training layers to learn a residual mapping $N(x) = H(x) - x$ rather than a direct mapping~$H(x)$~\cite{deep_res}. By utilizing shortcut connections that recast the output as $N(x) + x$, the network can easily favor identity mappings if the current representation is already optimal. This reparameterization ensures stable gradient flow and facilitates the training of very deep unrolled architectures, such as those required for dynamic CT reconstruction. An illustration of such a layer is shown in \cref{fig:res_layer}.

\section{Dataset} 

\noindent The development of robust neural network architectures for tomographic reconstruction requires diverse data that accounts for both the ideal reconstruction and the degradations inherent in real-world systems. 
While standard datasets like LoDoPaB-CT \cite{Leuschner2019TheLD} provide a baseline for quantum noise, they do not include the structural operator inexactness caused by mechanical instabilities. Consequently, we developed a specialized simulation framework to model dynamic perturbations $\Gamma$ in nano-CT environments. 

\subsection{Phantom Generation and Geometry}

\noindent To ensure, that the models learn generalized features rather than specific shapes, we generated a large-scale dataset of 2D phantoms $F \in \mathbb{R}^{M }$ with $M=255^2$. Each phantom is an assembly of a primary geometric structure and three subsidiary shapes. These shapes are randomly selected to be either ellipses or rectangles, with dimensions sampled from a uniform distribution:
\begin{itemize}
    \item Primary Shapes: Coverage of $4\%$ to $16.3\%$ (rectangles) or $3.2\%$ to $12.6\%$ (ellipses).
    \item Subsidiary Shapes: Coverage of $0.26\%$ to $0.96\%$ (rectangles) or $0.198\%$ to $0.78\%$ (ellipses).
\end{itemize}
The intensity values for these shapes are sampled from $\mathcal{U}(0.2, 0.8)$, while the background is set to zero. To increase spatial complexity, subsidiary shapes are placed within a larger radius (76 pixels) than the primary shape (26 pixels), occasionally resulting in shapes that are partially clipped by the primary shape. 

\subsection{Modeling Temporal Perturbations}

\noindent The core of our investigation involves the dynamic inverse problem, where the object moves relative to the scanner during the acquisition process. We ignore quantum noise ($\delta = 0$) to isolate the effects of the deformation $\Gamma$. We define two distinct motion models acting over the time interval $[0, T]$, discretized into $K=567$ time steps corresponding to the scan angles $\varphi_k \in [0, 180^\circ]$.

\subsubsection{Vibration Motion}

To simulate high-frequency mechanical instabilities, we employ a \emph{vibration} model ($d_V$) consisting of 38 overlapping damped sinusoidal waves. These $x$ and $y$ shifts are defined by functions 
$$\epsilon_x(t) = A_x \sin(2\pi\omega t) e^{-c t}, \quad \epsilon_y(t) = A_y \sin(2\pi\omega t) e^{-c t}.$$ 
The amplitudes $A$ follow a normal distribution $\mathcal{N}(7.8, 161.3)$, while frequency $\omega$ and damping $c$ are sampled uniformly. This shift is coupled with a minor random rotation $\xi_\gamma$ sampled from a truncated normal distribution. The total distortion is formalized as an affine mapping
$$d_V(t, (x', y')) = \left( OS_1 \cdot \mathrm{R}(\xi_\gamma, \epsilon_x, \epsilon_y) \cdot OS_2 \begin{pmatrix} x' \\ y' \\ 1 \end{pmatrix} \right)^T,$$
where $OS_1$ and $OS_2$ are translation matrices centering the rotation around the phantom's midpoint.
\subsubsection{Linear Shift Motion}
The \emph{linear} model ($d_L$) simulates a constant-velocity drift, common in thermal expansion or stage-drift scenarios. The displacement is linear over time and defined as 
$$\epsilon_x(t) = t \cdot \xi_x, \quad \epsilon_y(t) = t \cdot \xi_y.$$
The velocity $\xi$ is truncated such that the maximum displacement over the full scan duration does not exceed $40\%$ of the FOV. This results in a smooth, wave-like distortion in the sinogram domain, distinct from the jittered appearance of the vibration data.

\subsection{Data Generation and Partitioning}
\noindent The sinograms $g^\eta \in \mathbb{R}^{K \times L}$ are synthesized using the Operator Discretization Library (ODL) \cite{odl} and the ASTRA Toolbox \cite{astra_toolbox}. The implementation of this data generation requires a discrete formulation of the Radon transform and the forward projection. The field of view is discretized into pixels $S_m$, $m=1, ..., M$. 
For each beam $L_j:=L_{\omega (\varphi_k), s_l}$ the total length of the intersections of this beam with the pixels $S_m$ in the field of view is calculated by
$a_{j, m} = \text{Length of } (L_{j} \cap S_m)$, with $j=k\cdot l$.
For each scan angle $k$, the Radon transform $\R_d$ is applied to the specifically deformed state of the phantom $F^{\eta, t_k}\in\mathbb{R}^M$ at time instance $t_k$. 
The discrete (inexact) Radon transform is then be written as
\begin{equation}\label{eq:discrete_inexact_radon}
(\R^\eta_d f)(\omega(\varphi_k), s_l) = \sum_{m=1}^M  F^{\eta, t_k}_{m} a_{j, m},\qquad j=k\cdot l
\end{equation}
where $M=363$ represents the number of detector pixels. A more detailed explanation of the discrete inexact Radon transform can be found in \cite{Oberacker_2025}.
The dataset characteristics are summarized in Table \ref{tab:mydata_split}.

\begin{table}[h]
  \caption{Number of samples per data split for vibration and linear shift data.}
  \centering
  \begin{tabular}{lcccc}    
    \rowcolor{tablelighterblue} 
    object motion & Train & Val & Test & Total\\ \hline
    \hline
    Vibrations & 30\,489 & 1283 & 323 & 32\,095\\
    Linear shifts & 11\,631 & 369 & 323 & 12\,323 \\ \hline
  \end{tabular}
  \label{tab:mydata_split}
\end{table}

The linear shift dataset is smaller by design; because the motion is more deterministic (constant velocity), it is hypothesized that the neural network requires fewer samples to characterize the deformation than the highly stochastic vibration model.

\section{Optimization of Learned Methods}
\noindent This section outlines the evaluation of the proposed unrolled network architectures. The performance is assessed across two datasets: the vibration dataset (random motion) and the linear shift dataset (directional motion). Evaluation metrics include Peak Signal-to-Noise Ratio (PSNR) and Structural Similarity Index (SSIM), with a specific focus on the trade-offs between network depth and the risk of overfitting.


\subsection{Optimization of FLO Architecture on the Vibration Dataset}
\noindent The vibration dataset, characterized by complex random motion, serves as the primary benchmark for architectural tuning. The evaluation focuses on two sets of parameters: iteration block depth (denoted as \emph{depth sd} and \emph{depth $g^\eta$}) and SESOP-derived parameters (number of unrolled iterations $N$ and number of search directions $\# \text{sd}$).

\subsubsection{Depth and Overfitting Analysis}
Experiments on the FLO architecture $\mathcal{F}^\AA_\theta$ reveal that increased network depth correlates strongly with overfitting. To quantify this, we utilize a relative loss difference metric
$$\text{loss}_{\text{diff}} := \frac{|\text{train loss} - \text{val loss}|}{\text{train loss}}.$$
A value of $\text{loss}_{\text{diff}} < 0.5$ is established as the threshold for acceptable generalization. Detailed results of the initial experiments can be found in \cref{tab:sd2_it3_learned_op_vibrations}. Tests with $N=3$ and $\# \text{sd}=2$ demonstrate that while deeper configurations (e.g., depth $g^\eta=4$) may achieve superior PSNR (32.65 dB), they exhibit severe overfitting ($\text{loss}_{\text{diff}} = 3.18$). The optimal configuration for the FLO model is identified as depth $g^\eta=2$ and depth sd=1, achieving an SSIM of 0.921. This model maintains a negligible $\text{loss}_{\text{diff}}$ of 0.0018, demonstrating robust generalization.

\begin{table}[h]
\caption{Fully learned operator $\mathcal{N}^{\AA}_{\theta}(g^\eta)$ evaluated on validation dataset for measurements obtained from perturbed sinograms. Evaluated settings for depth of $g^\eta$ and search direction blocks with 3 unrolled iterations and 2 search directions. The best results are marked in bold, when the relative loss difference is too high the next best result is additionally underlined.}
\label{tab:sd2_it3_learned_op_vibrations}

\centering
\begin{tabular}{llccc}
    \rowcolor{tablelighterblue} 
    \multicolumn{5}{l}{N = 3, \# sd = 2} \\ 
    \rowcolor{tabledarkerblue} 
    &  &  &  & \\
    \rowcolor{tabledarkerblue} 
    \multirow{-2}{0.8cm}{depth $g^\eta$} & \multirow{-2}{0.8cm}{depth sd }& \multirow{-2}{*}{\textbf{PSNR [db]} }& \multirow{-2}{*}{\textbf{SSIM}}& \multirow{-2}{1cm}{$\text{loss}_{\text{diff}}$}\\
    \hline \hline 
1 & 1 & $32.36 \pm 2.73$ & $0.886 \pm 0.008$ & 0.3386 \\ 
  & 2 & $31.93 \pm 2.60$ & $0.836 \pm 0.005$ & 0.2128 \\ 
  & 3 & $17.41 \pm 3.00$ & $0.894 \pm 0.028$ & 0.5328 \\ 
  & 4 & \underline{$\mathbf{32.63 \pm 2.79}$} & $0.812 \pm 0.006$ & 0.3188 \\\hline 
2 & 1 & $31.84 \pm 2.54$ & $\mathbf{0.921 \pm 0.010}$ & 0.0018 \\ 
  & 2 & $32.43 \pm 2.62$ & $0.842 \pm 0.007$ & 0.3151 \\ 
  & 3 & $31.73 \pm 2.47$ & $0.813 \pm 0.008$ & 1.8827 \\ 
  & 4 & $31.27 \pm 2.60$ & $0.836 \pm 0.007$ & 0.2847 \\\hline 
3 & 1 & $19.29 \pm 0.50$ & $0.148 \pm 0.036$ & 4.7365 \\ 
  & 2 & $31.60 \pm 2.53$ & $0.857 \pm 0.006$ & 4.0593 \\ 
  & 3 & $30.91 \pm 2.49$ & $0.858 \pm 0.008$ & 0.3160 \\ 
  & 4 & $31.40 \pm 2.54$ & $0.806 \pm 0.008$ & 1.4375 \\\hline 
4 & 1 & $31.64 \pm 2.36$ & $0.675 \pm 0.048$ & 0.3207 \\ 
  & 2 & $32.56 \pm 2.77$ & $0.821 \pm 0.007$ & 0.3165 \\ 
  & 3 & $\mathbf{32.65 \pm 2.75}$ & $0.834 \pm 0.006$ & 3.1780 \\ 
  & 4 & $31.54 \pm 2.07$ & $0.622 \pm 0.010$ & 2.1048 \\ 
    \hline
\end{tabular}

\end{table}

\subsubsection{SESOP Parameter Tuning}
Analysis of the number of iterations versus search directions indicates that larger models do not inherently provide better reconstructions as can be seen in \cref{tab:learned_op_vibrations_sd_vs_iter}. For instance, a model with 4 search directions and 7 iterations achieved a PSNR of 32.06 dB but suffered from a massive $\text{loss}_{\text{diff}}$ of 7.95. The study concludes that the $N=3, \# \text{sd}=2$ configuration remains the most effective balance of performance and stability.
\begin{table}[h]
\caption{Fully learned operator $\mathcal{N}^{\AA}_{\theta}(g^\eta)$ evaluated on validation dataset for measurements obtained from perturbed sinograms. Evaluated settings for unrolled iterations and number of search directions with depth of $g^\eta$ of 2 and depth of search direction blocks of 1.  The best results are marked in bold, when the relative loss difference is too high the next best result is additionally underlined.}
\label{tab:learned_op_vibrations_sd_vs_iter}

\centering
\begin{tabular}{llccc}
    \rowcolor{tablelighterblue}
    \multicolumn{5}{l}{depth sd $=1$, depth $g^\eta=2$} \\ 
    \rowcolor{tabledarkerblue}
    & & & & \\
    \rowcolor{tabledarkerblue}
    \multirow{-2}{*}{\# sd} & \multirow{-2}{*}{N} & \multirow{-2}{*}{\textbf{PSNR [db]}} & \multirow{-2}{*}{\textbf{SSIM}} & \multirow{-2}{1cm}{$\text{loss}_{\text{diff}}$}\\ 
    \hline \hline 
2 & 3 & $31.84 \pm 2.54$ & $\mathbf{0.921 \pm 0.010}$  & $0.0018$\\ 
  & 4 & $30.44 \pm 2.15$ & $0.835 \pm 0.007$  & $0.2039$\\ 
  & 5 & $31.40 \pm 2.34$ & $0.784 \pm 0.011$  & $0.1962$\\\hline 
3 & 4 & $\mathbf{32.28 \pm 2.65}$ & $0.777 \pm 0.014$  & $4.1257$\\ 
  & 5 & $32.00 \pm 2.62$ & $0.861 \pm 0.006$  & $0.2752$\\ 
  & 6 & \underline{$\mathbf{32.06 \pm 2.71}$} & $0.809 \pm 0.007$  & $0.2198$\\\hline 
4 & 5 & $31.91 \pm 2.74$ & $0.816 \pm 0.007$  & $0.3327$\\ 
  & 6 & $31.99 \pm 2.54$ & $0.800 \pm 0.009$  & $15.2512$\\ 
  & 7 & $\mathbf{32.06 \pm 2.74}$ & $0.843 \pm 0.006$  & $7.9521$\\ 
    \hline
\end{tabular}

\end{table}


\subsection{Optimization of Radon Architecture on the Vibration Dataset}
\noindent The following summary evaluates the Radon Network ($\mathcal{F}_\theta^\mathcal{R}$), a physics-informed unrolled architecture that incorporates the Radon transform and its adjoint directly into the iterative refinement process.

\subsubsection{Depth and Overfitting Analysis}
The depth of the search direction blocks ($\textit{depth sd}$) was systematically evaluated to determine its impact on reconstruction quality and model stability. Unlike the FLO network, the Radon network does not utilize a separate $g^\eta$ block, focusing its learning capacity entirely on the search direction refinements. \\
Increasing the search direction depth to 2 was found to be highly beneficial for the Radon architecture. All experiments utilizing a depth of 2 achieved SSIM values exceeding 0.9, suggesting that a minimum level of non-linear complexity is required to effectively process the physics-based residuals generated by the Radon operator. While deeper models generally performed better, the evaluation monitored the relative loss difference between training and validation sets. It was observed that a configuration with a search direction depth of 2 offered the best trade-off, achieving superior Peak Signal-to-Noise Ratio (PSNR) values while maintaining a lower relative loss difference compared to significantly larger or deeper configurations.

\subsubsection{SESOP Parameter Tuning}
The Radon network’s performance is heavily influenced by parameters derived from the SESOP method: the number of unrolled iterations ($N$) and the number of search directions ($\#sd$). The architecture was tested with variations in iterations (ranging up to 4) and search directions (ranging up to 3). Given that the Radon transform is applied in every iteration, each additional step significantly increases the computational burden. The best-performing configuration was identified as 4 unrolled iterations and 3 search directions. This setup provided the network with sufficient historical search information to compensate for the dynamic motion artifacts reintroduced by the $\mathcal{R}^*\mathcal{R}$ operations. Increasing the number of search directions beyond this point led to diminishing returns and a prohibitive increase in training time without providing substantial improvements in image fidelity.

\begin{table}[h]
\caption{Quantitative evaluation of Radon network $\mathcal{F}^{\mathcal{R}}_{\theta}(g^\eta)$ evaluated on validation dataset for measurements obtained from perturbed sinograms. The best results are marked in bold, when the relative loss difference is too high the next best result is additionally underlined.}
\label{tab:radon_vibrations_sd_vs_iter}

\centering
\begin{tabular}{lllccc} 
    \rowcolor{tablelighterblue}
     &  &  &  &   & \\
    \rowcolor{tablelighterblue}
    \multirow{-2}{0.6cm}{depth sd} & \multirow{-2}{0.6cm}{\# sd}& \multirow{-2}{0.2cm}{N}& \multirow{-2}{*}{\textbf{PSNR [db]}}& \multirow{-2}{*}{\textbf{SSIM}}& \multirow{-2}{0.6cm}{$\text{loss}_{\text{diff}}$}\\
    \hline \hline 
1 & 2 & 3 & $30.41 \pm 2.77$ & $0.866 \pm 0.038$  & $0.0150$\\ 
  &   & 4 & $30.53 \pm 2.91$ & $0.869 \pm 0.037$  & $0.1019$\\ 
  & 3 & 4 & $30.57 \pm 2.85$ & $0.913 \pm 0.033$  & $0.1221$\\ 
  &   & 5 & $30.13 \pm 2.79$ & $0.864 \pm 0.025$  & $0.0474$\\\hline 
2 & 2 & 3 & $30.92 \pm 2.85$ & $0.949 \pm 0.015$  & $0.0758$\\ 
  &   & 4 & $31.33 \pm 2.77$ & $0.955 \pm 0.013$  & $0.0085$\\ 
  & 3 & 4 & $\mathbf{31.35 \pm 2.69}$ & \underline{$\mathbf{0.959 \pm 0.011}$}  & $0.0158$\\ 
  &   & 5 & $30.93 \pm 2.43$ & $\mathbf{0.961 \pm 0.011}$  & $0.1838$\\ 
    \hline
\end{tabular}

\end{table}

\subsection{Transferability of Configurations to Dataset with Different Motion}
\noindent To assess how well the developed architectures generalize across different types of dynamic noise, the optimal network settings identified in the vibration experiments were applied to models trained on the linear shift data. The FLO network utilized 3 unrolled iterations and 2 search directions (with a search direction depth of 1 and a $g^\eta$ depth of 2). Similarly, the Radon network was trained using 4 unrolled iterations and 3 search directions with a depth of 2. This approach offers insight into how architecturally robust the system remains when faced with various modeling inaccuracies arising from different types of motion.

\section{Evaluation of Learned Methods with Other Reconstruction Techniques}
\noindent The following section provides a comparative analysis of the proposed unrolled neural networks, the FLO Model and the Radon Network, against three established reconstruction techniques: Filtered Back Projection (FBP), Dremel method and SESOP-Kaczmarz method. This evaluation highlights the relative strengths of learned methods in mitigating motion-induced artifacts in dynamic CT data.

\subsection{Evaluation on Vibration Dataset}
\begin{figure}[h!]
  \centering
  \includegraphics[scale=0.5]{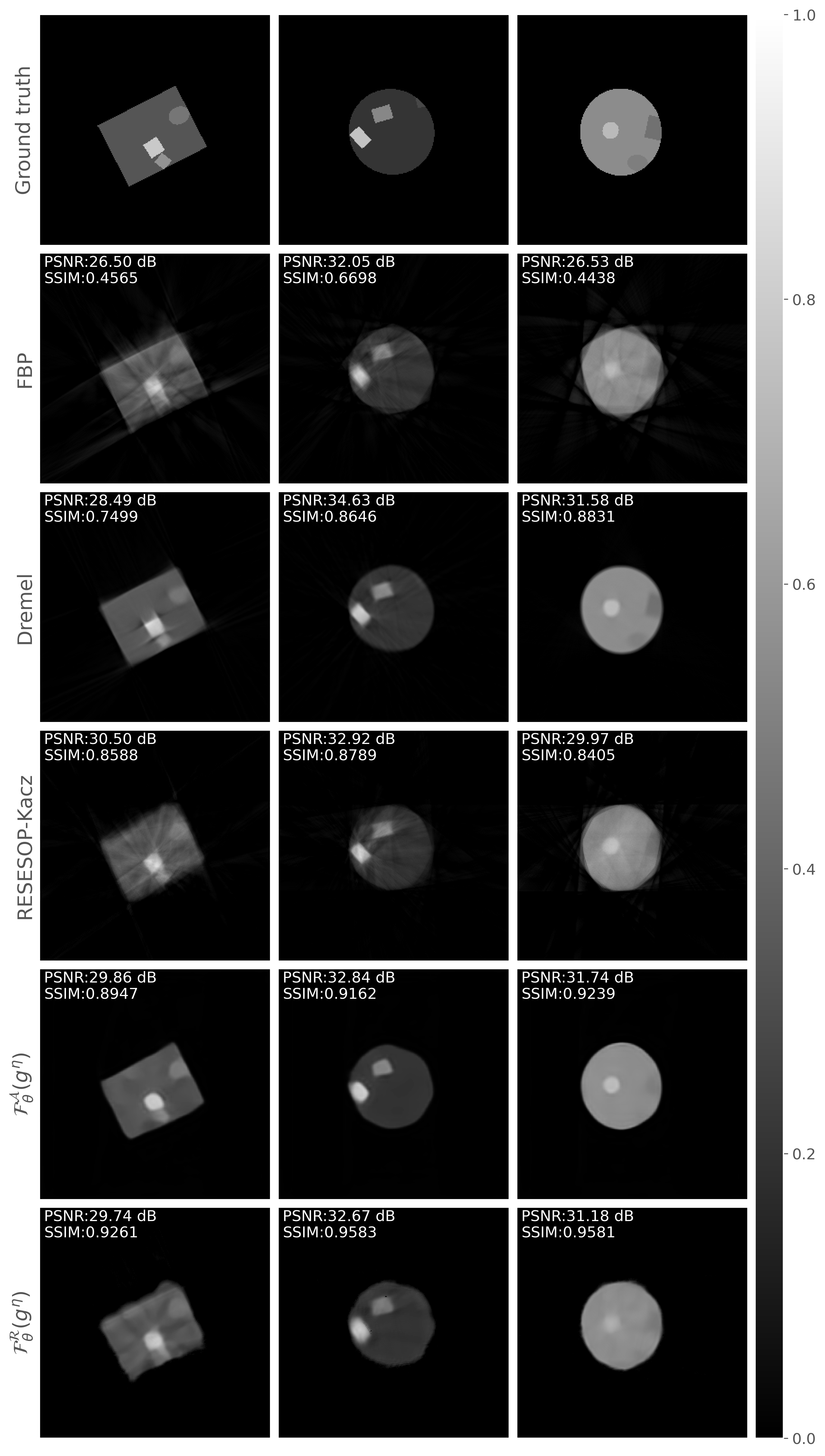}
  \caption{Reconstructions of measurement data perturbed with vibrations by all considered reconstruction methods.}
  \label{fig:exp1_all}
\end{figure}

\noindent A detailed examination of the reconstruction samples provides insight into how the FLO and Radon models handle dynamic noise. Visual inspection confirms that both models successfully reconstruct the target objects without generating the external streak artifacts typical of analytical methods.
Three samples from the vibration dataset were chosen to demonstrate the performance of all methods; see \cref{fig:exp1_all}.

The FLO model ($\mathcal{F}^\AA_\theta(g^\eta)$) demonstrates a superior ability to reproduce sharp edges and internal shapes. While it performs exceptionally well on elliptical structures, it struggles slightly with rectangular shapes, where corners appear less sharp than in the ground truth. In contrast, the Radon model ($\mathcal{F}^\mathcal{R}_\theta(g^\eta)$) produces less defined edges; rectangular corners are poorly reproduced, and elliptical boundaries often exhibit a "wave-like" distortion. 

\begin{figure}[h!]
  \centering
  \includegraphics[scale=0.5]{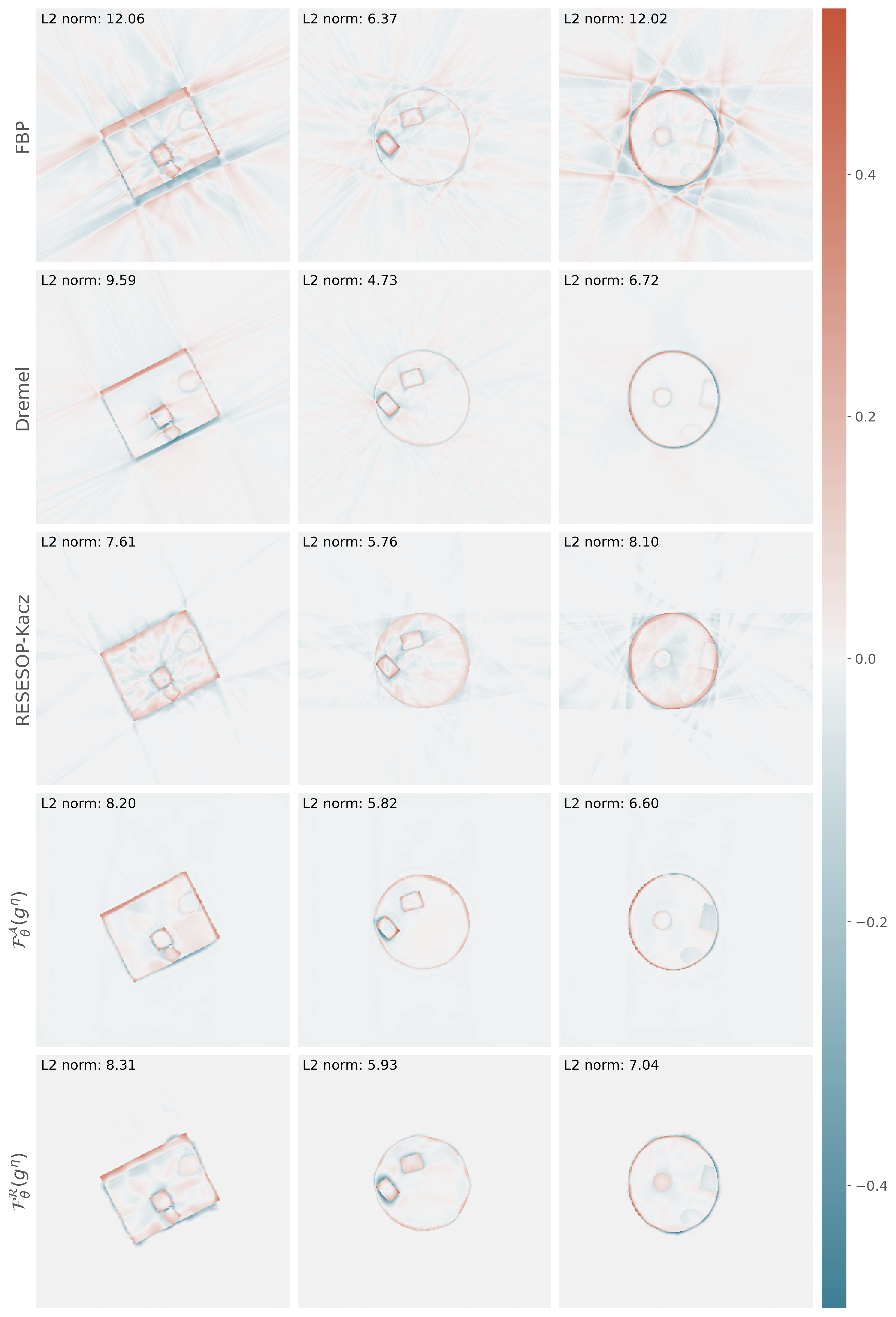}
  \caption{Differences between ground truths and reconstructions of measurement data perturbed with vibrations by all considered reconstruction methods.}
  \label{fig:exp1_all_diff}
\end{figure}

\Cref{fig:exp1_all_diff} shows the differences between the ground truths and the reconstructions, where red indicates an overestimated density and blue indicates an underestimated density.
Analysis of these difference maps reveals distinct error patterns. The FLO model tends to shift the object slightly, as evidenced by the blue and red edges on opposite sides of the structure. The Radon model, however, shows mixed colors along the edges, indicating geometric distortion rather than a simple shift. Interestingly, while the Radon model achieved slightly higher SSIM values, the FLO model secured slightly higher PSNR scores and produced visually cleaner and less blurry results. 

To contextualize the performance of the neural networks, they were compared against non-learned methods. The quantitative results, summarized in \cref{tab:exp1_all}, indicate that the learned architectures achieve significantly higher fidelity across both PSNR and SSIM metrics.
\begin{table}[h]
\caption{Quantitative evaluation on vibration test datasets for measurements obtained from perturbed sinograms. The best results for both measures are highlighted.}
  \label{tab:exp1_all}
\centering
\begin{tabular}{lcc}
    \rowcolor{tablelighterblue} 
    & \textbf{PSNR [dB]} & \textbf{SSIM} \\
    \hline \hline
    FBP & $27.93 \pm 2.82$ & $0.511 \pm 0.104$ \\ 
    Dremel & $30.99 \pm 3.01$ & $0.829 \pm 0.076$ \\ 
    RESESOP-Kacz. with $\eta$ & $30.67 \pm 2.74$ & $0.865 \pm 0.102$ \\ 
    Fully learned op. net & $\mathbf{31.61 \pm 2.54}$ & $0.920 \pm 0.010$ \\ 
    Radon net & $31.00 \pm 2.61$ & $\mathbf{0.958 \pm 0.011}$ \\ 
    \hline
\end{tabular}

\end{table}
Among the non-learned techniques, the Dremel method is the most effective producing fewer streak artifacts than FBP or SESOP-Kaczmarz. However, its artifacts are still visible outside the object, and its edges remain significantly more blurred than those produced by the FLO model. FBP and SESOP-Kaczmarz suffer from substantial streaking both inside and outside the object, making them less suitable for precision dynamic imaging.

Ultimately, the FLO model provides the most balanced performance, combining well-defined edges and accurate density mapping without the computational artifacts seen in the Radon network or the severe streaking characteristic of traditional methods.

\newpage

\newpage

\newpage

\begin{figure}[h!]
  \centering
  \includegraphics[scale=0.5]{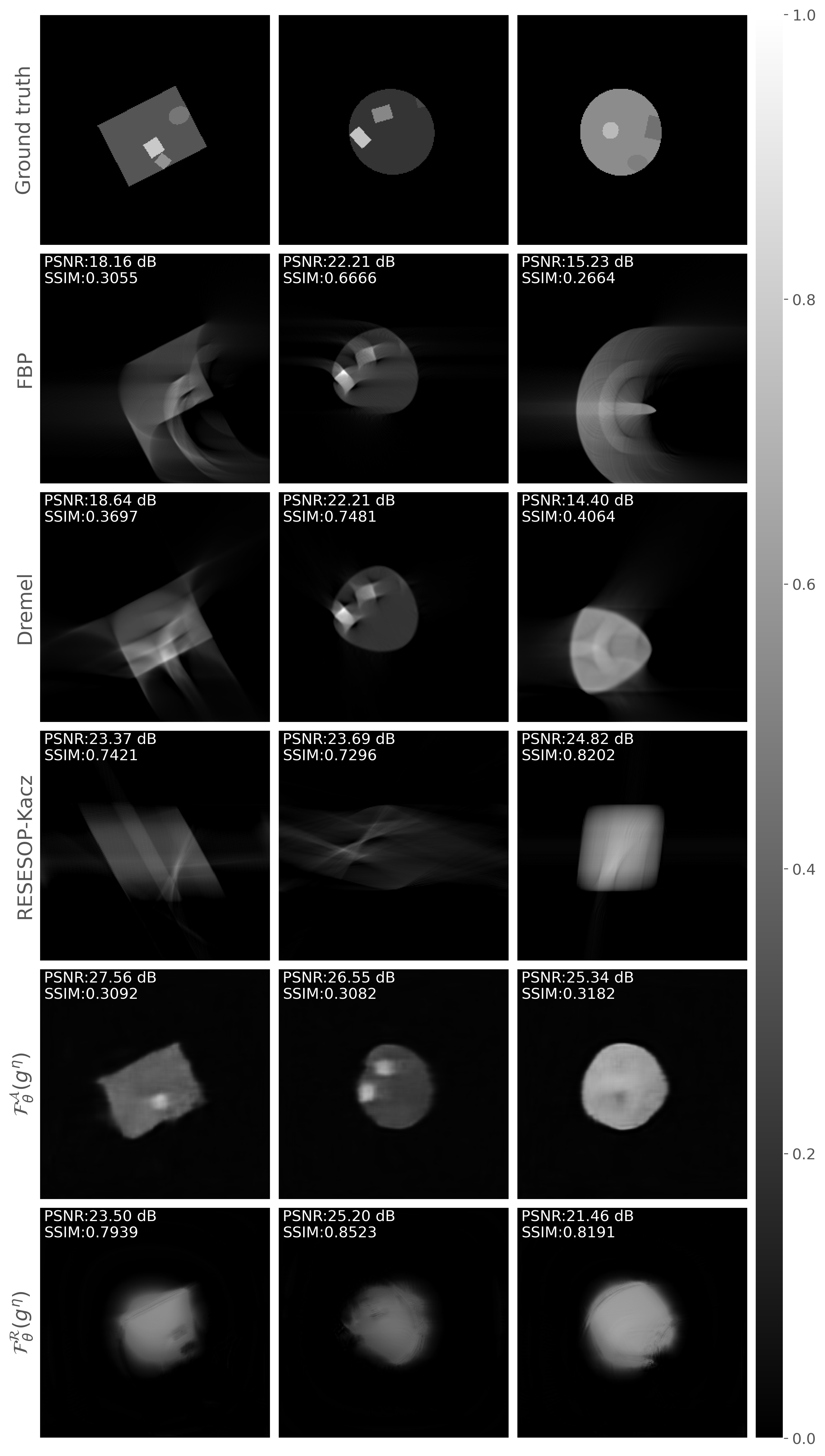}
  \caption{Reconstructions of measurement data perturbed with linear shifts by all considered reconstruction methods.}
  \label{fig:exp2_all}
\end{figure}

\newpage

\begin{figure}[H]
  \centering
  \includegraphics[scale=0.5]{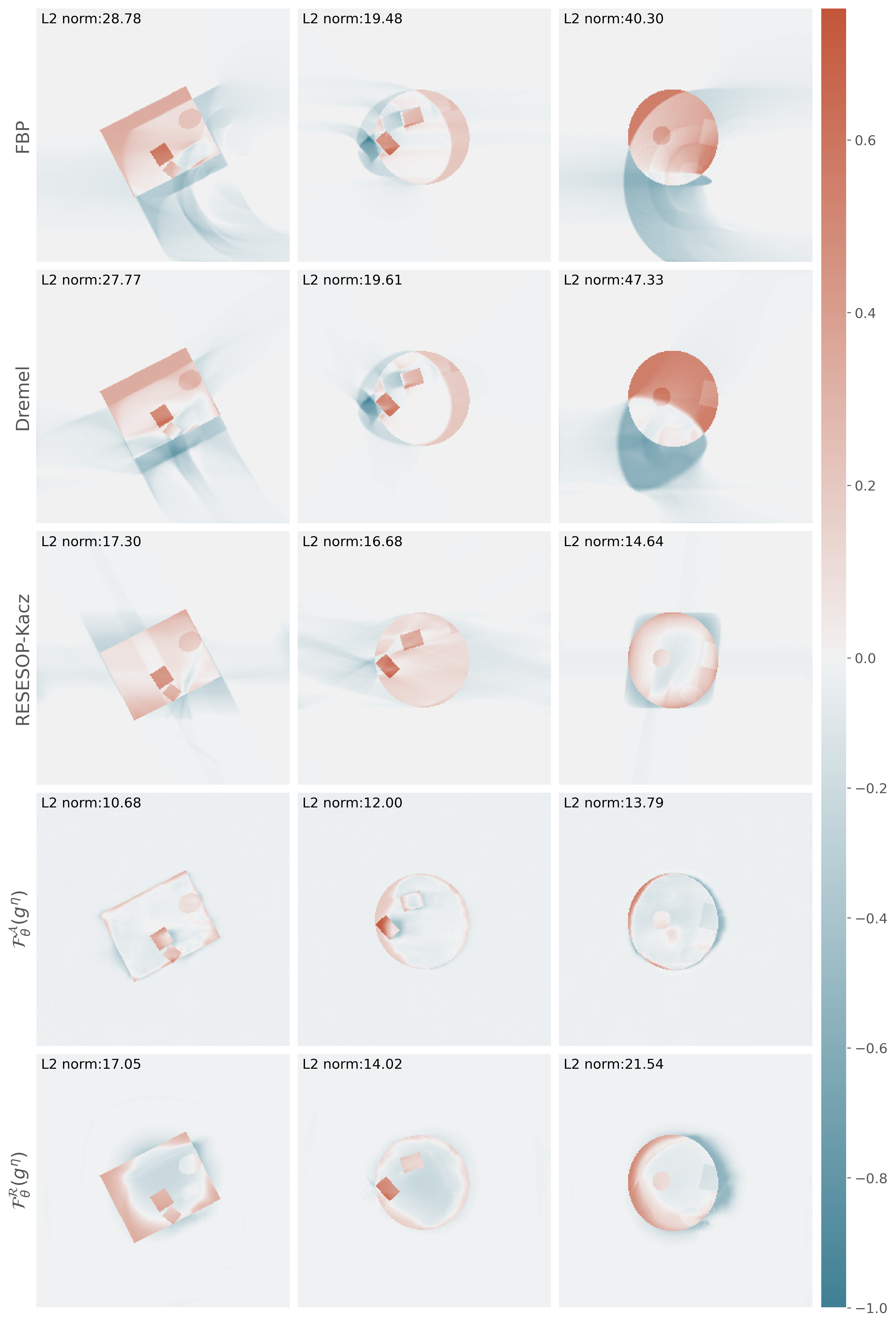}
  \caption{Differences between ground truths and reconstructions of measurement data perturbed with linear shifts by all considered reconstruction methods.}
  \label{fig:exp2_all_diff}
\end{figure}

\newpage

\subsection{Evaluation on Linear Shift Dataset}
\noindent The comparison of reconstruction models under linear shift perturbations reveals significant challenges in recovering high-quality images, with no single model producing entirely satisfactory results.

Both learned methods use an initial reconstruction as input for $f_0$. Due to its fast computation FBP reconstructions are used. In spite of the extreme motion in the linear shift dataset, the FBP reconstructions provided an initialization that improved PSNR by several decibels compared to random starts. More details on this analysis can be found in \cite{Oberacker_2025}.
The quantitative results of the reconstruction techniques on this dataset can be found in \cref{tab:exp2_all}.
\begin{table}[h]
\caption{Quantitative evaluation on linear shift test dataset for measurements obtained from perturbed sinograms. The best results for both measures are highlighted.}
  \label{tab:exp2_all}
\centering
\begin{tabular}{lcc}
    \rowcolor{tablelighterblue} 
    & \textbf{PSNR [dB]} & \textbf{SSIM} \\
    \hline \hline
    FBP & $17.71 \pm 3.63$ & $0.395 \pm 0.144$ \\ 
    FBP shifted & $20.46 \pm 4.22$ & $0.480 \pm 0.118$ \\ 
    Dremel & $17.41 \pm 3.79$ & $0.482 \pm 0.150$ \\ 
    Dremel shifted& $21.50 \pm 4.84$ & $0.581 \pm 0.132$ \\ 
    RESESOP-Kacz. with $\eta$ & $22.88 \pm 3.37$ & $0.777 \pm 0.055$ \\ 
    Fully learned op. net & $\mathbf{24.95 \pm 2.41}$ & $0.653 \pm 0.016$ \\ 
    Radon net & $23.03 \pm 2.69$ & $\mathbf{0.819 \pm 0.017}$ \\ 

    \hline
\end{tabular}

\end{table}
The reconstructions of three samples from the linear shift dataset are illustrated in \cref{fig:exp2_all}.

The Radon model performed poorly from a visual standpoint, typically producing circular artifacts with only fragmented edges visible. Despite these visual shortcomings, the Radon model surprisingly achieved the highest SSIM values. Analysis of difference maps confirms that these reconstructions are visually inferior to those of the FLO model, suggesting that in this specific context, SSIM may not accurately reflect perceptual quality or density accuracy. The FLO model achieved the highest PSNR values among all methods (24.95 dB).

Among the non-learned techniques, the Dremel method provides the most visually acceptable results in scenarios with low motion; however, it quickly succumbs to smearing and distortion as the shift intensity increases. Analytical FBP exhibits even more pronounced artifacts, while the SESOP-Kaczmarz method fails entirely to resolve discernible shapes, proving incapable of compensating for the systematic linear drift. 

Evaluation shows that the objects in non-learned reconstructions are frequently reconstructed for different time points. Ground truths were always generated for the first time point. While an automated image translation via pixel cross-correlation was applied to evaluate object quality independent of location, the distortions were often too severe for effective correction.
Such shifts are clearly visible in the difference maps in \cref{fig:exp2_all_diff}.

Ultimately, the FLO model emerges as the superior architecture for handling extreme modeling inexactness. Even in cases where the linear shift is large enough to render the objects unrecognizable to traditional methods, the FLO model consistently reconstructs the primary phantom shapes and recovers some of the internal features. Despite the occasional presence of large-scale artifacts, its ability to maintain density accuracy and structural integrity under severe perturbations makes it the most robust solution for dynamic CT reconstruction in linear shift scenarios.

\section{Conclusion}

We developed data-driven reconstruction methods for dynamic CT imaging at a nanoscale (nano-CT). The methods rely on neural networks that are trained in a supervised manner and whose architectures are unrolled iterative regularization schemes. We specifically suggested a fully learned operator (FLO) model, meaning that the unrolled neural network has no information about the exact (dynamic) forward model at all, as well as a Radon network, where we used the (static) Radon transform as an approximation to the exact model. The performances of the FLO model and the Radon network were compared with other methods like the classical Filtered Backprojection (FBP), the Dremel method, which has been developed specifically for dynamic CT, and iterative schemes such as SESOP and SESOP-Kaczmarz. We used simultaed data, where the motion models were gicen by vibration and linear shift models. For both motion models the learned methods proved superior to the classical schemes with the FLO showing the best balanced results. For the linear shift motion none of the methods yielded convincing results except the FLO model, where the objects are still clearly visible what proves the FLO method to be a robust solution approach for dynamic CT reconstruction, even at a non scale.

\bibliographystyle{abbrv}
\bibliography{references}


\bigskip
\bigskip
\bigskip

\appendix

%

\section{Training Details}\label{app:training_details}

\noindent The training configuration for the unrolled networks was designed to balance convergence speed with model stability and generalization. The settings are summarized as follows.\\

\underline{Loss function:} 
The mean squared error (MSE) was selected as the loss function. Since CT data (both sinograms and phantoms) are inherently bounded between $0$ and $1$, the risk of outliers is minimal, allowing MSE to provide faster convergence than the mean absolute error. \\

\underline{Optimization:} Optimization was conducted using the AdamW algorithm, provided by PyTorch \cite{pytorch}, which decouples weight decay from the gradient update for more effective regularization. The hyperparameters were set to default values: learning rate ($\epsilon$): $0.001$, decay rates ($\rho_1, \rho_2$): $0.9$ and $0.999$, mini-batch size: $32$ (constrained by GPU memory limits).\\

\underline{Regularization and Convergence:}
To enhance generalizability and prevent overfitting, the training pipeline incorporated three distinct regularization methods:\\[1ex]

\begin{enumerate}
    \item Batch Normalization: Integrated directly into the architecture to stabilize internal activations. \\[1ex]
    \item $L^2$ Regularization: Controlled via the AdamW weight decay parameter, set to $\lambda = 0.01$. \\[1ex]
    \item Early Stopping: Employed for extended training runs (beyond 30 epochs). The process was terminated if the validation loss failed to improve by more than $\Delta = 0.00015$ over a patience of 5 epochs. \\
\end{enumerate}
This configuration ensured that the high-dimensional inverse problem was solved efficiently while maintaining the structural integrity of the reconstructed objects.

\section{SESOP-Kaczmarz Settings}
The SESOP-Kaczmarz algorithm for this research is implemented in Python, relying on two primary parameters: the operator inexactness $\eta$ and the measurement noise $\delta$. In this specific experimental setup, the operator inexactness is treated as fully known, representing an optimal but non-realistic scenario. To simplify the modeling, the algorithm utilizes the maximum norm of $\eta$ per scan angle rather than point-wise values, which slightly relaxes the regularization process. This inexactness parameter serves two critical roles: it limits the size of the reconstruction stripes $H_n^{\eta, \delta, \rho}$ and informs the discrepancy principle used for terminating the iteration.

The discrepancy principle is further tuned using a safety factor $\tau$ set to 1.00001 to ensure the reconstruction error can effectively fall below the established noise level. Although the synthetic measurement data lacks actual quantum noise, a baseline $\delta$ value of 0.001 is introduced to artificially provide a noise level for regularization purposes. The forward and adjoint operators, representing the Radon transform and its backprojection, are managed via the ASTRA toolbox \cite{astra_toolbox}.

Final termination of the algorithm is governed by both the discrepancy principle and a fixed iteration limit. For the evaluations conducted, the maximum number of iterations ($N$) is fixed at 20, which frequently results in the process terminating before the discrepancy principle is formally satisfied to maintain consistency across the comparative study.

\section{Dremel Settings}
The implementation of the Dremel method for this study was developed in collaboration with researchers from the University of Bremen \cite{NN_inexactOp}, with all parameters and procedures optimized specifically for the vibration datasets. A key technical configuration is the use of a super-resolution factor of 2, which up-samples the projections to a size of 726 points. While the original literature suggested using an edge filter for pre-processing low-contrast signals, this study found that a combination of mean subtraction and zero padding yielded superior results. To handle the output, the post-processing phase involves cropping the correlation signal back to its original dimensions and applying a specific weighting to compensate for the bias toward smaller shifts introduced by the zero padding.

The operational parameters for the algorithm were set with a learning rate and relaxation parameter both equal to 1, meaning no relaxation was applied. Notably, the operator correction is initiated from the very first iteration, diverging from the standard approach of including a warm-up phase without correction.  Finally, instead of using a discrepancy principle to decide when to stop, the algorithm is configured to terminate after a fixed limit of 32 iterations, a value determined by assessing performance on the validation data.

\section{Filtered Back Projection Settings}
The FBP implementation for this research is facilitated by the ASTRA toolbox \cite{astra_toolbox}, employing the default Ram-Lak filter.


\section{Architecture Details} \label{app:arch_details}

The following two pictures illustrate architecture details for the FON and the Radon network.

\begin{figure}[H]
  \centering
  \includegraphics[width=3in]{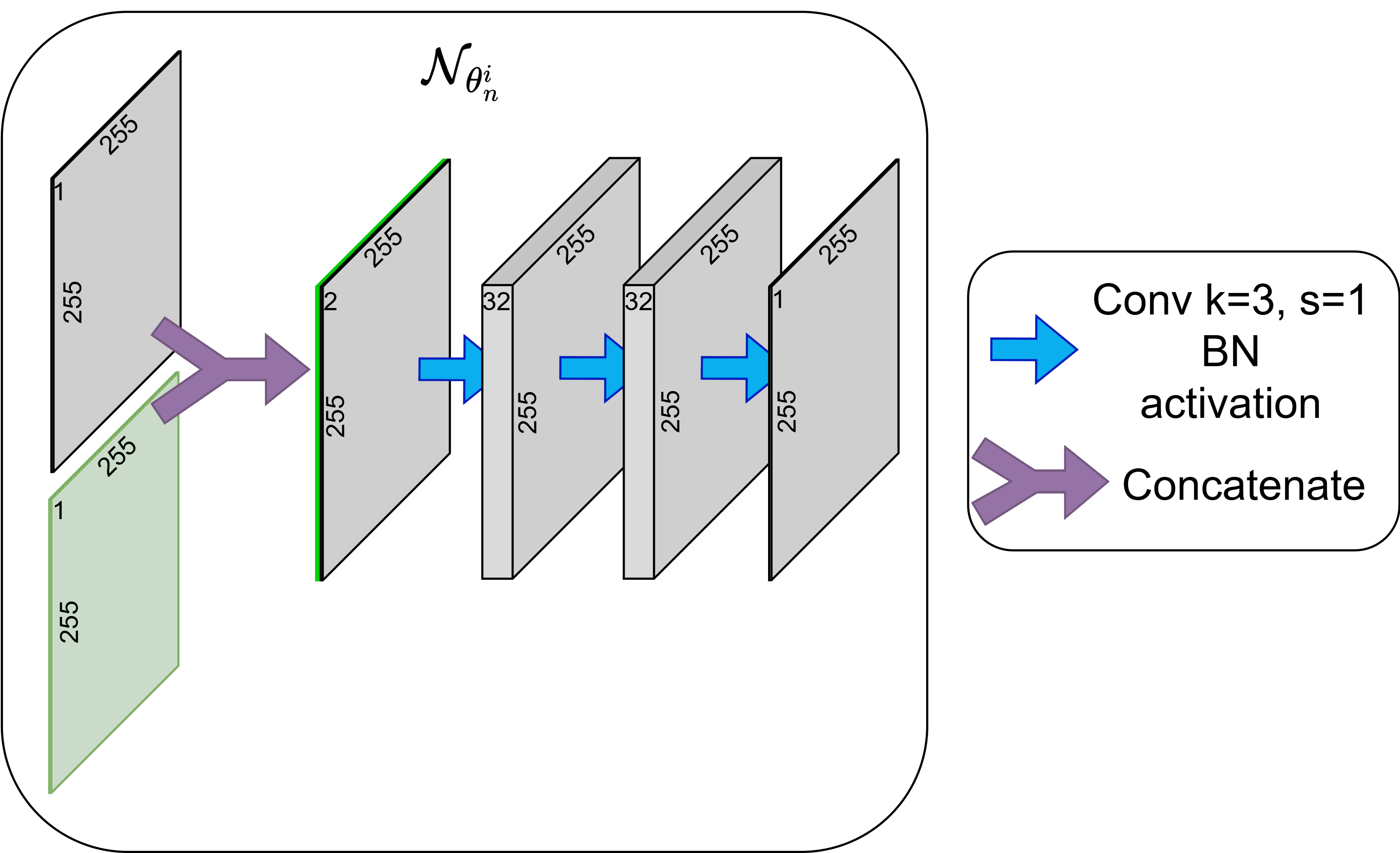}
  \caption{This block uses the Radon transform of reconstruction $f_i$, shown in green, and $f_i$ itself, shown in gray, to feed into several convolutional layers combined with batch normalization and activation. This block has a depth of 2, as 2 intermediate layers with 32 channels are used. The output data always has a channel depth of 1.}
  \label{app:adjoint_sesop_iteration}
\end{figure}

\begin{figure}[H]
  \centering
  \subfloat[]{
  \includegraphics[width=5in]{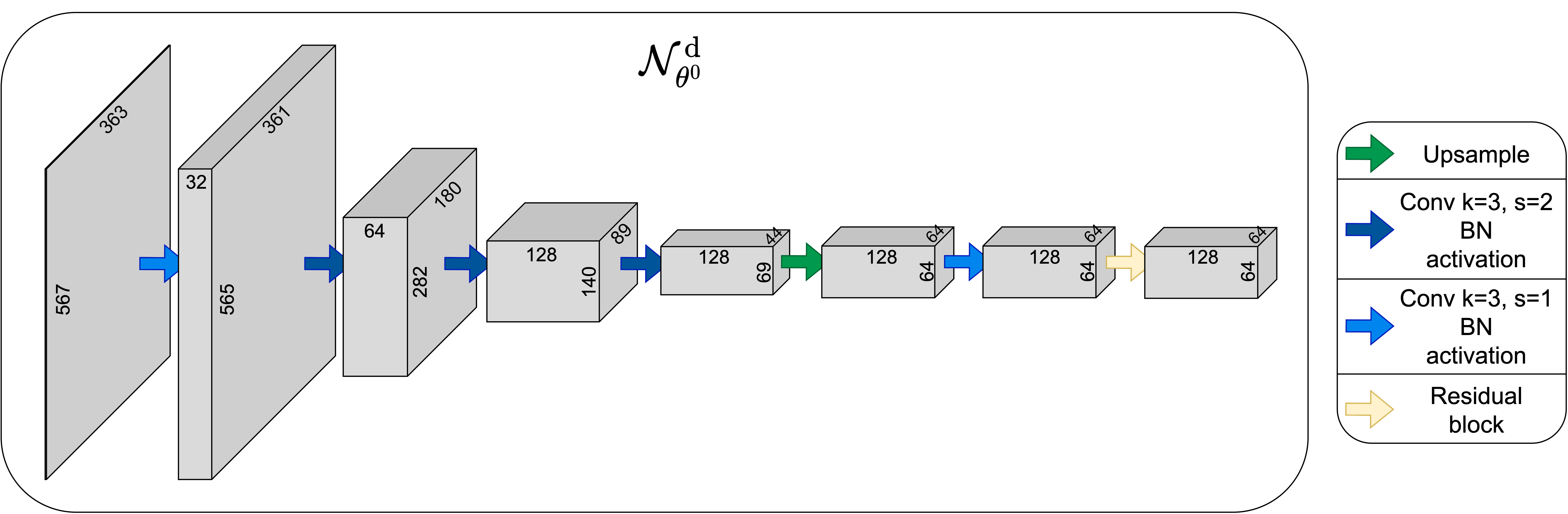}%
  \label{app:N_0(g)}}
\\
\subfloat[]{
  \includegraphics[width=5in]{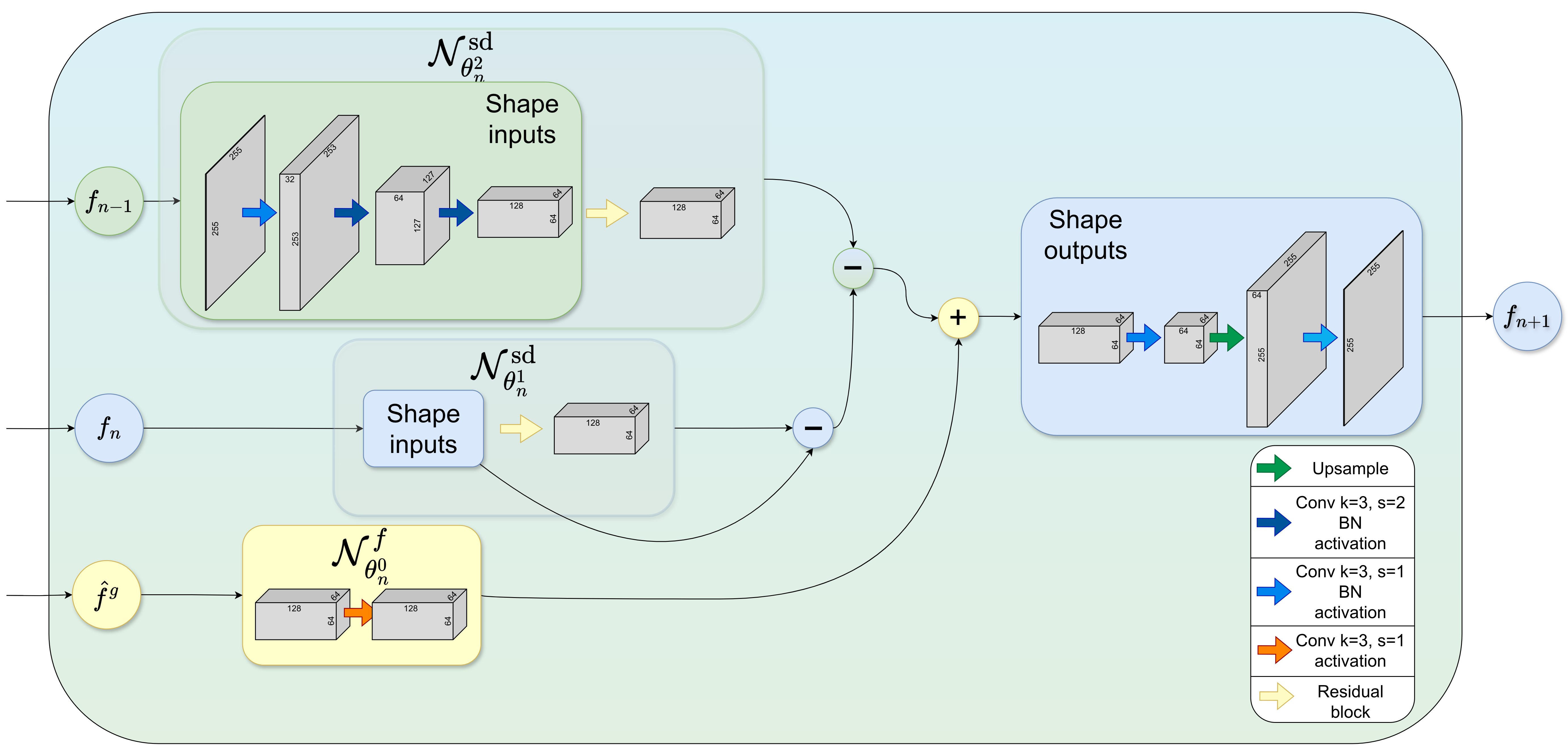}%
  \label{app:detail_sesop_block}}
  \caption{(a) The measurement data undergoes several convolutional transformations. Blue arrows represent different CBA layer configurations. The green arrow symbolizes an upsampling step. The yellow arrow references a residual block. (b) A more detailed overview of the internal architecture of the SESOP-block. The search direction data is first reshaped and then passed into residual blocks. The output of the data processing block $\mathcal{N}^{\text{d}}_{\theta^0}$ is further transformed in $\mathcal{N}^{f}_{\theta^0_n}$. After performing the operations given by the SESOP iteration the output is transformed into its required shape.}
\end{figure}

\end{document}